\documentclass{birkmult}
\usepackage{amssymb}
\usepackage{amsmath}
\usepackage{amsthm}

\theoremstyle{plain}

\newtheorem{theorem}{Theorem}[section]
\newtheorem{lemma}[theorem]{Lemma}
\newtheorem{fact}[theorem]{Fact}
\newtheorem{proposition}[theorem]{Proposition}
\newtheorem{corollary}[theorem]{Corollary}
\newtheorem{observation}[theorem]{Observation}
\newtheorem{example}[theorem]{Example}

\theoremstyle{definition}

\newtheorem{definition}[theorem]{Definition}
\newtheorem{remark}[theorem]{Remark}
\newtheorem{problem}[theorem]{Problem}

\newcommand{\R}{\mathbb{R}}
\newcommand{\N}{\mathbb{N}}

\newcommand{\E}{\mathbb{E}}
\newcommand{\p}{\mathbb{P}}

\def\epsilon{\varepsilon}
\def\rho{\varrho}
\def\phi{\varphi}
\def\theta{\vartheta}

\begin{document}
%-------------------------------------------------------------------------
% editorial commands: to be inserted by the editorial office
%
%\firstpage{1}
%\volume{228}
%\DOI{003-0001}
%\seriesextraline{This is the Concrete Title of this Book\br H.E. R and S.T.C. W, Eds.}
%\firstpage{1}
\title[Delta-semidefinite and delta-convex quadratic forms]
 {Delta-semidefinite and delta-convex\\ quadratic forms in Banach spaces}

%----------Author 1
\author[N.~Kalton]{Nigel Kalton}

\address{%
         Department of Mathematics\\
         University of Missouri-Columbia\\
         Columbia, MO 65211, U.S.A.}

\email{nigel@math.missouri.edu}

%----------Author 2
\author[S.V.~Konyagin]{Sergei V. Konyagin}

\address{%
         Department of Mechanics and Mathematics\\
         Moscow State University\\
         Moscow, 119992,
         Russia}

\email{konyagin@ok.ru}

%----------Author 3
\author[L.~Vesel\'y]{Libor Vesel\'y}
\address{%
         Dipartimento di Matematica\\
         Universit\`a degli Studi di Milano\\
         Via C.~Saldini, 50\\
         20133 Milano,
         Italy}
\email{vesely@mat.unimi.it}

\thanks{The first author was supported by NSF grant DMS-0555670.  The second author was supported by the Russian Foundation for
Basic Research, Grant~05-01-00066, and by Grant NSh-5813.2006.1.
The third author was supported in part by the Ministero
dell'Universit\`a e della Ricerca of Italy.}



\date{May 15, 2006}

\begin{abstract}
A continuous quadratic form (``quadratic form'', in short) on a
Banach space $X$ is: (a) {\em delta-semidefinite} (i.e., representable as
a difference of two nonnegative quadratic forms) if and only if the
corresponding symmetric linear operator $T\colon X\to X^*$ factors through a Hilbert space;
(b) {\em delta-convex} (i.e., representable as a difference of two continuous convex
functions) if and only if $T$ is a UMD-operator.
It follows, for
instance, that each quadratic form on an infinite-dimensional
$L_p(\mu)$ space ($1\le p \le \infty$) is: (a) delta-semidefinite if{f} $p \ge 2$;
(b) delta-convex if{f} $p>1$. Some other related results concerning
delta-convexity are proved and some open problems are stated.
\end{abstract}

%%% ----------------------------------------------------------------------
\maketitle
%%% ----------------------------------------------------------------------
%\tableofcontents

%%%%%%%%%%%%%%%%%%%%%%%%%%%%%%%%%%%%%%%%%%%%%%%%%%%%%%%%%%%%%%%%%%%

\section*{Introduction}

Let $X$ be a real Banach space. Recall that
a function $q\colon X\to \R$ is
a {\em continuous quadratic
form} (more precise would be ``continuous purely quadratic form'')
if there exists a
continuous bilinear form $b\colon X\times X\to\R$ such that $q(x)=b(x,x)$
for each $x\in X$.

In the present paper, we are interested mainly in the following two
isomorphic properties of $X$.

\begin{enumerate}
\item[(D)] {\em Each continuous quadratic form on $X$ is {\em delta-semidefinite},
i.e., it can be represented as a difference of two nonnegative
continuous quadratic forms.}
\item[(dc)] {\em Each continuous quadratic form on $X$ is {\em delta-convex},
i.e., it can be represented as a difference of two continuous convex functions.}
\end{enumerate}
Since nonnegative quadratic forms are convex, (D) always
implies (dc). The reverse implication is not true, as we shall see in
Section~\ref{S:three}.

In Section~\ref{S:one}, we characterize delta-semidefinite continuous quadratic forms
on $X$ as
precisely those whose corresponding symmetric linear operator
$T\colon X\to X^*$ is factorizable through a Hilbert space. This leads, via
known results on factorizability, to sufficient conditions for a
Banach space $X$ to satisfy (D). The characterization also implies that
the property (D) passes to quotients, and
the spaces $\ell_p$, $1\le p<2$, do not satisfy (D).

In Section~\ref{S:two}, we use $X$-valued Walsh-Paley martingales to prove
that a continuous quadratic form on $X$ is delta-convex if and only if
the corresponding symmetric linear operator is a UMD-operator. It
follows that $\ell_1$ not only fails (D) but it also fails (dc).

In Section~\ref{S:three}, we discuss relationships between the
properties (D), (dc) and the following property.

\smallskip
\leftline{(Cdc)\ \ {\em Each $C^{1,1}$ function on $X$ is delta-convex.}}
\smallskip

\noindent
It is easy to see that also (Cdc) implies (dc). We show that (dc) and (Cdc) pass 
to quotients.
For each of the
properties (D), (dc), (Cdc), we characterize those $p$'s in $[1,\infty]$ for
which an infinite-dimensional $L_p(\mu)$ space satisfies the property
(Theorem~\ref{T:L_p}).
It follows that (dc) implies neither (D) nor (Cdc). (The latter should be compared with
a result from \cite{Du-Ve-Za} which says that all Banach space-valued
quadratic mappings on $X$ are delta-convex if and only if all Banach space-valued
$C^{1,1}$ mappings on $X$ are delta-convex.)
We also solve a
problem from \cite{Ve-Za} by proving existence of a
function $f$ whose compositions with all ``delta-convex curves''
(in the sense of \cite{Ve-Za}) are
delta-convex while $f$ is not locally delta-convex. Some of these
counterexamples use a result by M.~Zelen\'y \cite{Zel}. Finally, we show
that the property (dc) is not stable with respect to direct sums, and we
state some open problems.

%%%%%%%%%%%%%%%%%%%%%%%%%%%%%%%%%%%%

\section{Delta-semidefinite quadratic forms}\label{S:one}

In what follows, the term ``operator'' means ``bounded linear
operator''.
Recall that an operator $T\colon X\to X^*$ is called {\em symmetric} if
$\langle Tx,y\rangle = \langle Ty,x\rangle$ for all $x,y\in X$
(equivalently: $T^*=T$ on $X$).

It is easy to see that the formula
\begin{equation}\label{E:qT}
q(x)=\langle Tx,x\rangle
\end{equation}
defines a one-to-one correspondence between the continuous quadratic
forms $q$ on $X$ and the symmetric operators $T\colon X\to X^*$.\\
(Indeed, if $q$ is generated by a continuous bilinear form $b$, it is
generated also by the symmetric bilinear form $\frac{b(x,y)+b(y,x)}{2}$.
Moreover, there is a unique symmetric bilinear form $b$ that generates
$q$; this follows from the formula
\begin{equation}\label{E:polar}
2b(x,y)=b(x+y,x+y)-b(x,x)-b(y,y)
\end{equation}
valid for symmetric $b$. The rest follows from the well-known
one-to-one
correspondence, via the formula
$b(x,y)=\langle
Tx,y\rangle$,
between the continuous bilinear forms $b$ on $X\times X$ and
the operators $T\colon X\to X^*$.)\\
If (\ref{E:qT}) holds for each $x\in X$, we say that {\em $T$ generates $q$}.

The formula (\ref{E:polar}) also implies the following

\begin{fact}\label{F:p'}
Each continuous quadratic form $q$ on $X$ is everywhere Fr\'echet differentiable.
Moreover, its Fr\'echet derivative at $x$ is given by $q'(x)=2Tx$
where $T$ is the symmetric operator
that generates $q$.
\end{fact}

The following theorem characterizes delta-semidefinite continuous quadratic forms.
(Recall that a continuous quadratic form is called {\em delta-semidefinite}
if it is the difference of two nonnegative continuous quadratic forms.)
An operator $T\colon X\to Y$ is said to be {\em factorizable
through $Z$} if there exist operators
$A\colon X\to Z$ and $B\colon Z\to Y$ such that $T=BA$.

\begin{theorem}\label{T:kalton}
Let $q$ be a continuous quadratic form on a Banach space $X$, and
$T\colon X\to X^*$ be the symmetric operator that generates $q$. Then
the following assertions are equivalent:
\begin{enumerate}
\item[(i)] $q$ is delta-semidefinite;
\item[(ii)] there exists a continuous quadratic form $p$ on $X$, such
that $|q|\le p$;
\item[(iii)] $T$ is factorizable through a Hilbert space.
\end{enumerate}
\end{theorem}

\begin{proof}\ \\
$(iii)\Rightarrow(i)$. If $T=BA$ where $A\colon X\to H$ and
$B\colon H\to X^*$ are operators, and $H$ is a Hilbert space, then we have
\[
q(x)=\langle BAx,x\rangle=\langle Ax, B^*x\rangle_H=
{\textstyle\frac{1}{4}}\|Ax+B^*x\|^2_H-
{\textstyle\frac{1}{4}}\|Ax-B^*x\|^2_H
\]
which shows that $q$ is
difference of two nonnegative quadratic forms.

\smallskip\noindent
$(i)\Rightarrow(ii)$. If $q=p_1-p_2$ where $p_i$ ($i=1,2$) are
nonnegative continuous quadratic forms, then
$|q|\le p_1+p_2=:p$.

\smallskip\noindent
$(ii)\Rightarrow(iii)$.
Let (ii) hold, and let $S\colon X\to X^*$ be the symmetric operator such
that $p(x)=\langle Sx,x\rangle$. The function
\[
[\cdot,\cdot]\colon X/\mathrm{Ker}(S) \times X/\mathrm{Ker}(S) \to\R,
\quad
[\xi,\eta]:=\langle Sx,y\rangle\
\text{where $x\in\xi$, $y\in\eta$,}
\]
is well-defined, bilinear, symmetric, and $[\xi,\xi]\ge0$ for each
$\xi\in X/\mathrm{Ker}(S)$. Moreover, if $p(x)=0$ for some $x\in X$,
then $x$ is a minimizer for $p$, and hence $0=p'(x)=2Sx$ by
Fact~\ref{F:p'}. In other words, $[\xi,\xi]=0$ implies $\xi=0$.
Consequently,
$[\cdot,\cdot]$ is an inner product on $X/\mathrm{Ker}(S)$.
Let $H$ be the completion of the inner product space
$\bigl(X/\mathrm{Ker}(S),[\cdot,\cdot]\bigr)$.

Consider the operator $J=iQ\colon X\to H$ where
$Q\colon X\to X/\mathrm{Ker}(S)$ is the quotient map and
$i\colon X/\mathrm{Ker}(S)\hookrightarrow H$ is the inclusion map.
($J$ is continuous since
$\|Qx\|^2_H=\langle Sx,x\rangle\le\|S\|\cdot\|x\|^2$ for all $x\in X$.)

If $x\in\mathrm{Ker}(S)$, then $p(x)=q(x)=0$. Since $p+q$ is a
nonnegative quadratic form generated by the symmetric operator $T+S$,
the same argument as above shows that $Tx+Sx=0$. This proves that
$\mathrm{Ker}(S)\subset\mathrm{Ker}(T)$. Consequently, the operator
\[
T_0\colon X/\mathrm{Ker}(S)\to X^*,\quad T_0\xi:=Tx\
\text{where $x\in\xi$},
\]
is well-defined. We claim that $T_0$ is continuous also in the norm
generated by the inner product $[\cdot,\cdot]$. To prove this, consider
$\xi\in X/\mathrm{Ker}(S)$ and $y\in X$ such that
$\|\xi\|_H \le1$ and $\|y\| \le1$. Fix $x\in\xi$, and denote $\eta=Qy$.
Then
$\|\eta\|^2_H=\langle Sy,y\rangle\le\|S\|$, and
$
%\begin{align}
\bigl|\langle T_0\xi,y\rangle\bigr| \le
\bigl|\langle Tx,y\rangle\bigr| =
\frac{1}{2}\bigl|q(x+y)-q(x)-q(y)\bigr|\le
\frac{1}{2}\bigl[p(x+y)+p(x)+p(y)\bigr]=
\frac{1}{2}\bigl[\|\xi+\eta\|^2_H +\|\xi\|^2_H+ \|\eta\|^2_H\bigr]\le
\frac{1}{2}\bigl[(1+\|S\|^{1/2})^2+1+\|S\|\bigr].
%\end{align}
$

Thus $T_0$ has a unique extension to an operator from $H$ into $X^*$;
let us denote it by $T_0$ again. Then $T=T_0J$ is the desired
factorization through $H$.
\end{proof}

\begin{corollary}\label{C:quotient}
\begin{enumerate}
\item[(a)] A Banach space $X$ has the property (D) (see Introduction) if and only
if each symmetric operator $T\colon X\to X^*$ is factorizable through a
Hilbert space.
\item[(b)] The property (D) passes to quotients, and hence also to
complemented subspaces.
\end{enumerate}
\end{corollary}

\begin{proof}
(a) follows immediately from Theorem~\ref{T:kalton}.
Let us show (b). Let $X$ satisfy (D), and let $L$ be a closed
subspace of $X$. Let $T\colon X/L\to (X/L)^*=L^\perp$ be a symmetric operator.
Consider the operator $S=iTQ\colon X\to X^*$ where $Q\colon X\to X/L$ is
the quotient map, and $i\colon L^\perp\to X^*$ is the inclusion
isometry. Since $Q^*=i$, we have
$\langle Sx,y\rangle=\langle T(Qx),(Qy)\rangle$ ($x,y\in X$), which shows that
$S$ is symmetric. By (a), $S$ is factorizable through a Hilbert space.
Now, Proposition~7.3 in \cite{DJT} implies that $T$ factors through a
Hilbert space, too.
\end{proof}

Operators that are factorizable through a Hilbert space were intensively
studied (the main reference is \cite{Pisier}, see also \cite{DJT}), and there exist
many sufficient conditions for factorizability of all operators
between two given spaces.
Thus, by Theorem~\ref{T:kalton}, we obtain various
sufficient conditions for
validity of the property (D) (defined in Introduction); we collect them in the
following theorem.

For the classical notion of modulus of smoothness, see e.g. \cite{Lin-Tz}.
For the notion of type and cotype, see e.g.\ \cite{Lin-Tz} or \cite{DJT}.
We shall need the following notion of second order differentiability,
studied in \cite{Bor-Noll}.

\begin{definition}
\begin{enumerate}
\item[(a)]
Let $f$ be a continuous convex function on a Banach space $X$. We say
that $f$ is {\em second order differentiable}\
at a point $x_0\in X$ if
there exist $x_0^*\in X^*$ and a continuous quadratic form $q$ on $X$
such that, for each $v\in X$,
\[
f(x_0+tv)=f(x_0)+x^*_0(v)t+q(v)t^2+o(t^2)\qquad \hbox{as\ } t\to0.
\]
\item[(b)]
A {\em second order differentiable norm} is a norm which is second order
differentiable at each nonzero point.
\end{enumerate}
\end{definition}

\begin{remark}
It follows from results in \cite{Bor-Noll} that a norm $\|\cdot\|$ on
$X$ is second order differentiable if{f} it is Fr\'echet (equivalently:
G\^ateaux) smooth and its derivative
$\|\cdot\|'\colon X\setminus\{0\}\to X^*$ is ${\rm weak}^*$-G\^ateaux
differentiable.
\end{remark}

\begin{theorem}\label{T:D}
Let $X$ be a Banach space.
Each continuous quadratic form on $X$ is delta-semidefinite (and hence
delta-convex), provided {\em at least one} of the following conditions is
satisfied.
\begin{enumerate}
\item[(a)] $X$ has type 2.
\item[(b)] $X^*$ has cotype 2, and $X$ has the approximation property.
\item[(c)] $X^*$ has cotype 2, and $X$ does not contain $\ell_1(n)$'s
      uniformly.
\item[(d)] $X^*$ has cotype 2, and $X$ is a Banach lattice.
\item[(e)] $X=C(K)$ for some compact space $K$.
\item[(f)] $X=L_p(\mu)$ for $2\le p\le\infty$ and some positive measure
      $\mu$.
\item[(g)] $X=c_0(I)$ for some set $I$.
\item[(h)] $X$ admits a uniformly smooth renorming with modulus of
      smoothness of power type 2
      (i.e., $\rho_X(\tau)\le a\tau^2$ for some $a>0$).
\item[(i)] $X$ has the Radon--Nikod\'ym property and admits an
       equivalent second order differentiable norm.
\end{enumerate}
\end{theorem}

\begin{proof}
(a) follows from Corollary~3.6 and Proposition~3.2 in \cite{Pisier}.
\\
(b): see Theorem~4.1 in \cite{Pisier}.
\\
(c) follows from (a) by \cite[Corollary~2.5]{Pi1}.
\\
(d) follows from Theorems 8.17 and 8.11 in \cite{Pisier}.
\\
(e) follows e.g.\ from (d) since $C(K)^*$ has cotype 2 (see~\cite[p.34]{Pisier}).
\\
(f): the case $p<\infty$ follows from (a) (see~\cite[p.73]{Lin-Tz});
     the case $p=\infty$ follows from (e) (see~\cite[Theorem~1.b.6]{Lin-Tz}).
\\
(g) follows from (e) by Corollary~\ref{C:quotient}(b), since
$c_0(\Gamma)$ is a closed hyperplane in $c(\Gamma)=C(K)$ where $K$ is
the one-point compactification of the discrete set $\Gamma$.
\\
(h) follows from (a) (see\ Theorem~1.e.16 in \cite{Lin-Tz}).
\\
(i) follows from (h) by the following reasoning.
If the norm of $X$ is second order differentiable,
then this norm is Lipschitz-smooth at each point of $S_X$ by
\cite{Bor-Noll}; this implies (by \cite[Lemma~2.4]{FWZ})
that the gradient of the norm is pointwise Lipschitz at each point
of $S_X$. By  \cite[Corollary~III.2]{DGZ}, if $X$ has also the
RNP, then it satisfies (h).
\end{proof}

Results in \cite[Section$\,$5]{SSTT} imply that, 
for each $1<p<2$,  there exists an operator
$U\colon\ell_p\to\ell_{p^*}$ (where $\frac{1}{p}+\frac{1}{p^*}=1$) such that $U$ is
not factorizable through a
Hilbert space.
Since $U$, constructed in \cite{SSTT}, is also symmetric, we obtain one
more corollary to Theorem~\ref{T:kalton}.

\begin{corollary}\label{C:badp}
The space $\ell_p$ fails the property (D) whenever $1\le p<2$.\\
{\rm
(The case of $p=1$ follows from Corollary~\ref{C:quotient}(b) and from
the well-known fact that each separable Banach space is isometric to a
quotient of $\ell_1$.)
}
\end{corollary}

%%%%%%%%%%%%%%%%%%%%%%%%%%%%%%%%%%%%

\section{Delta-convex quadratic forms}\label{S:two}

Let $X$ be a Banach space. Recall that a continuous function $\phi\colon X\to\R$ is
{\em delta-convex} if it is the difference of two continuous convex
functions on $X$. It is easy to see that $\phi$ is delta-convex if and
only if there exists a (necessarily convex) continuous function $\psi$ on $X$
such that both $\pm\phi +\psi$ are convex. Every such function
$\psi$ is called a {\em control function} for $\phi$. Denoting
$$
\Delta^2\phi(x,y):=\phi(x+y)+\phi(x-y)-2\phi(x)\,,\qquad
x,y\in X,
$$
it is easy to see that $\psi$ is a control function for $\phi$ if and
only if $|\Delta^2\phi(x,y)| \le \Delta^2\psi(x,y)$ for all $x,y\in X$.

Since every nonnegative quadratic form is convex, each
delta-semidefinite quadratic form is delta-convex. As we shall see in
Section~3, the converse is not true in general.

In this section, we use $X$-valued Walsh-Paley martingales to study delta-convexity
of quadratic forms. We recall all needed definitions and properties to
make our exposition self-contained.

\medskip

Let $n\ge1$ be an integer, $\Gamma=\{-1,1\}$, $f\colon\Gamma^n\to X$.
Then the {\em expectation} of $f$ is defined as
$\E f=2^{-n}\sum_{\eta\in\Gamma^n} f(\eta)=\int_{\Gamma^n} f\,d\p$,
where $\p=\p_n$ is the uniformly distributed probability measure on
$\Gamma^n$.

For $0\le k\le n$, consider the $\sigma$-algebra
$\Sigma_k=\{A\times \Gamma^{n-k}: A\subset\Gamma^k\}.$
Obviously, a function $f\colon\Gamma^n\to X$ is $\Sigma_k$-measurable if
and only if $f$ depends only on the first $k$ coordinates (in particular,
all $\Sigma_0$-measurable functions are constant). For this reason, we
sometimes
view $\Sigma_k$-measurable functions on $\Gamma^n$ as functions on
$\Gamma^k$.

For $f\colon\Gamma^n\to X$ and $0\le k\le n$, the {\em conditional expectation
of $F$ w.r.t.\ $\Sigma_k$} is the $\Sigma_k$-measurable function
$\E(f|\Sigma_k)\colon \Gamma^n\to X$ which has the same integral (w.r.t.\ $\p$)
as $f$ over each element of $\Sigma_k$. It is easy to see that it is given by
\[\textstyle
\E(f| \Sigma_k)(\omega)=\int_{\Gamma^{n-k}}f(\omega,\cdot)\,d\p_{n-k}
\,,\qquad \omega\in\Gamma^k.
\]
Note that $\E(f|\Sigma_0)\equiv\E f$, $\E(f|\Sigma_n)=f$, and
$\E(\E(f|\Sigma_k))=\E f$.

In this paper, we consider only Walsh-Paley martingales of finite length.

\begin{definition}
Let $X$ be a Banach space. An {\em $X$-valued Walsh-Paley martingale}
is any
finite sequence $(f_0,\ldots, f_n)$ of $X$-valued functions on
$\Gamma^n$ such that
$f_k=\E(f_n|\Sigma_k)$ for each $0\le k\le n\,;$
or equivalently, each $f_k$ is $\Sigma_k$-measurable, and
$$\textstyle
f_k(\omega)=\frac{1}{2}f_{k+1}(\omega,-1)+\frac{1}{2}f_{k+1}(\omega,1)
\quad\text{whenever $0\le k<n$ and $\omega\in\Gamma^k$.}
$$
Given a Walsh-Paley martingale $(f_0,\ldots,f_n)$, the corresponding
{\em martingale differences} are the functions
$df_k=f_k-f_{k-1}$\ \  $(1\le k\le n).$
\end{definition}

\begin{remark}\label{O:mart}
Let $(f_0,\ldots,f_n)$ be an $X$-valued Walsh-Paley martingale. The above definition
easily implies the following properties.
\begin{enumerate}
\item[(a)] $\E (df_k | \Sigma_j)=0$ whenever $0\le j< k\le n$.
\item[(b)] $(f_0,\ldots,f_n,f_n,\ldots,f_n)$ is a Walsh-Paley martingale
no matter how many times $f_n$ is repeated.
\item[(c)] $(Tf_0,\ldots,Tf_n)$ is a $Y$-valued Walsh-Paley martingale whenever
$T\colon X\to Y$ is linear.
\item[(d)] If $0\le m< n$ and $\overline{\omega}\in\Gamma^m$, then the
finite
sequence $(g_0,\ldots,g_{n-m})$, where
$g_k:=f_{m+k}(\overline{\omega},\cdot)\colon\Gamma^{n-m}\to X$, is a Walsh-Paley martingale.
\item[(e)] $df_k(\omega)=\frac{1}{2}\omega_k
\bigl[f_k(\omega_1,\ldots,\omega_{k-1},1)-f_k(\omega_1,\ldots,\omega_{k-1},-1)\bigr]$
whenever $1\le k\le n$ and $\omega\in\Gamma^k$.
\item[(f)] $f_{k-1}(\omega)\pm df_k(\omega)=f_k(\omega_1,\ldots,\omega_{k-1},\pm\omega_k)$
whenever $1\le k\le n$ and $\omega\in\Gamma^k$.
\end{enumerate}
\end{remark}

A finite sequence $(\epsilon_1,\ldots,\epsilon_n)$ of functions on $\Gamma^n$
is said to be {\em predictable} if $\epsilon_k$ is $\Sigma_{k-1}$-measurable
for each $1\le k\le n$.

\begin{lemma}\label{L:mart}
Let $(f_0,\ldots,f_n)$ be an $X$-valued Walsh-Paley martingale.
\begin{enumerate}
\item[(a)] If $\phi\colon X\to\R$ is a function, then
$\E\sum_{k=1}^n \Delta^2\phi(f_{k-1},df_k)=2\E\phi(f_n)-
2\phi(f_0)$.
\item[(b)] If $1\le k\le n$ and $w\colon\Gamma^n\to\R$ is $\Sigma_{k-1}$-measurable, then
$\E (w\,df_k)=0$.
\item[(c)] If $(g_0,\ldots,g_n)$ is an $X^*$-valued Walsh-Paley martingale and
$(\epsilon_1,\ldots,\epsilon_n)$ is a predictable sequence of real-valued functions, then
$\E\langle df_k,\epsilon_j dg_j\rangle=0$ whenever $k\ne j$ and $k,j\in\{1,\ldots,n\}$.
\item[(d)] $(\E \max\limits_{0\le k\le n}\|f_k\|^p)^{1/p}\le\frac{p}{p-1}(\E\|f_n\|^p)^{1/p}$
for each real $p>1$.
\end{enumerate}
\end{lemma}

\begin{proof}
(a) Recall that $f_0$ is constant.
By Remark~\ref{O:mart}(f), the left-hand side equals
\begin{multline*}\textstyle
\sum_{k=1}^n \int_{\Gamma^k} \bigl[
\phi(f_k(\omega))+\phi(f_k(\omega_1,\ldots,\omega_{k-1},-\omega_k))-2\phi(\omega)
\bigr]\,d\p_k(\omega)
\\ \textstyle
=2\sum_{k=1}^n \bigl( \E\phi(f_k) - \E\phi(f_{k-1})\bigr)=
2\E\phi(f_n) - 2\E\phi(f_0).
\end{multline*}
(b) follows easily from Remark~\ref{O:mart}(a).
\\
(c) Let, e.g., $k<j$. Obviously, $\langle df_k,\epsilon_jdg_j\rangle$ is
$\Sigma_{j}$-measurable. For each fixed $\omega\in\Gamma^{j-1}$, we have
$$\E (\langle df_k,\epsilon_jdg_j\rangle|\Sigma_{j-1})(\omega)=
\E \langle df_k(\omega),\epsilon_j(\omega) dg_j(\omega,\cdot)\rangle
=\langle df_k(\omega), \epsilon_j(\omega) \E dg_j(\omega,\cdot)\rangle=0$$
by
Remark~\ref{O:mart}(d),(a). Thus
$\E \langle df_k,\epsilon_j dg_j\rangle=
\E\bigl(\E (\langle df_k,\epsilon_j dg_j\rangle|\Sigma_k)\bigr)=0$. The case
$k>j$ is similar.
\\
(d) It is easy to verify that $h(\eta):=\tilde{h}(\eta) v(\eta)$ works,
where $\tilde{h}\colon\Gamma^n\to[0,\infty)$ satisfies
$\E\tilde{h}^2=1$ and $(\E\|g\|^2)^{1/2}=\E(\tilde{h}\cdot\|g\|)$, and
$v(\eta)\in S_X$ is such that $\|g(\eta)\|\le 2\langle
v(\eta),g(\eta)\rangle$.
\\
(e) For real-valued martingales, this is the well-known
Doob's $L_p$-inequality (see e.g.\ \cite[Theorem~14.11]{williams}). In the general case,
consider the Walsh-Paley martingale $(g_0,\ldots,g_n)$ given by
$g_k=\E(\|f_n\|\,|\Sigma_k)$. Then
$\|f_k\|=\| \E(f_n|\Sigma_k)\|\le g_k$ and $\|f_n\|=g_n$.
Hence, using the scalar case, we get
$(\E \max\limits_{0\le k\le n}\|f_k\|^p)^{1/p}\le
(\E \max\limits_{0\le k\le n}g_k^p)^{1/p}\le
\frac{p}{p-1} (\E g_n^p)^{1/p}=
\frac{p}{p-1}(\E \|f_n\|^p)^{1/p}$.
\end{proof}

\begin{lemma}\label{1}
Let $\phi$ be a continuous real function on a Banach space $X$.
In order that $\varphi$ is delta-convex it is necessary and
sufficient that there is a continuous function
$\rho\colon X\to [0,\infty)$
such that if $(f_0,\ldots,f_n)$ is an $X$-valued Walsh-Paley martingale then
\begin{equation}\label{E:1}
\mathbb E \sum_{k=1}^n |\Delta^2\varphi(f_{k-1},df_k)|\le \mathbb E\rho(f_n).
\end{equation}
Moreover, in this case, $\rho$ can always be taken convex.
\end{lemma}

\begin{proof}
Let $\phi$ be delta-convex with a control function $\psi$.
By adding a suitable affine function, we can (and do) suppose that
$\psi\ge0$.
Now Lemma~\ref{L:mart}(a) implies
\begin{align*}
\textstyle
\mathbb E \sum\limits_{k=1}^n |\Delta^2\varphi(f_k,df_k)|&
\textstyle
\le \mathbb E \sum\limits_{k=1}^n \Delta^2\psi(f_k,df_k)\\
&= 2\mathbb E\psi(f_n)-2\psi(f_0)\\
&\le 2\mathbb E\psi(f_n).
\end{align*}
Thus \eqref{E:1} holds with $\rho=2\psi$ (which is convex).

Conversely, if \eqref{E:1} holds, we may define
\begin{equation}\label{control}
\textstyle
\psi(x)=
\frac12\,\inf\Big\{ \mathbb E\rho(f_n)-
\mathbb E \sum\limits_{k=1}^n |\Delta^2\phi(f_{k-1},df_k)|\Big\}
\end{equation}
where the infimum is taken over all $X$-valued Walsh-Paley martingales with $f_0= x.$

Suppose $x,u\in X$, $y=x+u,z=x-u$ and $\epsilon>0.$
Using Remark~\ref{O:mart}(b) and the definition of $\psi$,
pick $X$-valued Walsh-Paley martingales $(f_0,\ldots,f_n)$ and $(g_0,\ldots,g_n)$
such that $f_0= y$, $g_0= z$ and
\begin{align*}
\psi(y)&\textstyle
>\frac12\,\mathbb E\rho(f_n)-\frac12\,\mathbb E \sum\limits_{k=1}^n
|\Delta^2\phi(f_{k-1},df_k)|-\epsilon\\
\psi(z)&\textstyle
>\frac12\,\mathbb E\rho(g_n)-\frac12\,\mathbb E \sum\limits_{k=1}^n
|\Delta^2\phi(g_{k-1},dg_k)|-\epsilon.
\end{align*}
Form a new Walsh-Paley martingale $(h_0,\ldots,h_{n+1})$ by setting $h_0=x$
and, for $1\le k\le n+1$,
$$ h_k(\eta_1,\ldots,\eta_{n+1})=\begin{cases} f_{k-1}(\eta_2,\ldots,\eta_{n+1}) \qquad \text{if } \eta_1=1\\
g_{k-1}(\eta_2,\ldots,\eta_{n+1})\qquad \text{if } \eta_1=-1.\end{cases}$$
Note that, for example,
$$\textstyle
\E \rho(h_{n+1})=
\int\limits_{\{\eta_1=1\}} \rho(h_{n+1})\,d\p_{n+1}+\int\limits_{\{\eta_1=-1\}} \rho(h_{n+1})\,d\p_{n+1}
=\frac12\E \rho(f_n) +\frac12\E \rho(g_n).
$$
Thus
\begin{align*}
2\psi(x)
&\textstyle
\le \E\rho(h_{n+1})- \E \sum\limits_{k=1}^{n+1} |\Delta^2\phi(h_{k-1},dh_{k})|\\
&\textstyle
=\textstyle{\frac12}\,\E\rho (f_n) + {\frac12}\,\E\rho (g_n) \\
&\qquad\quad\textstyle
- |\Delta^2\phi(x,u)|-
{\frac12}\,\E \sum\limits_{k=1}^n |\Delta^2\phi(f_{k-1},df_k)|-
{\frac12}\,\E \sum\limits_{k=1}^n |\Delta^2\phi(g_{k-1},dg_k)| \\
&\le \psi(x) + \psi(y) + 2\epsilon - |\Delta^2\phi(x,u)|.
\end{align*}
Since $\epsilon>0$ was arbitrary, it follows that
$|\Delta^2\phi(x,u)|\le \Delta^2\psi(x,u)$ whenever $x,u\in X$.
Thus $\psi$ is a midconvex (or Jensen convex) function which is locally bounded since
$0\le\psi\le\rho/2$. Consequently (see \cite[p.215]{Ro-Va}), $\psi$ is a
continuous convex function. Thus $\psi$ is a control function for
$\phi$.
\end{proof}

\begin{observation}\label{obs}
For each $g\colon\Gamma^n\to[0,\infty)$ and $p>0$, we have
$$\sum_{j=1}^\infty 2^{jp}\,\p(g>2^{j}) \le \frac{2^p}{2^p-1}\,\E g^p.$$
\rm Indeed,
$
\int g^p\,d\p
\ge \sum\limits_{j=1}^\infty \int\limits_{\{2^{j-1}<g\le 2^j\}}g^p\,d\p
\ge \sum\limits_{j=1}^\infty 2^{(j-1)p}\,\bigl[\p(g>2^{j-1})-\p(g>2^j)\bigr]
= \sum\limits_{j=0}^\infty 2^{jp}\, \p(g>2^j) - \sum\limits_{j=1}^\infty 2^{(j-1)p} \,\p(g>2^j)
\ge (1-2^{-p}) \sum\limits_{j=1}^\infty 2^{jp}\, \p(g>2^j).
$
\end{observation}

\bigskip

Let $p>0$. Recall that a function $\phi\colon X\to \R$ is called {\em positively $p$-homogeneous}
if
$\phi(tx)=t^p\phi(x)$ whenever $t\ge0$, $x\in X$.

\begin{lemma}\label{2}
Suppose $p>1$.  A continuous positively
$p$-homogeneous function $\phi\colon X\to\R$ is
delta-convex if and only if there is a constant
$C$ such that for all $X$-valued Walsh-Paley martingales $(f_0,\ldots,f_n)$ we have
\begin{equation}\label{E:2}
\mathbb E\sum_{k=1}^n|\Delta^2\varphi(f_{k-1},df_k)|\le C\mathbb
E\|f_n\|^p.
\end{equation}
\end{lemma}

\begin{proof}
Assume $\phi$ is delta-convex. Let $\rho$ be the corresponding continuous
function from Lemma~\ref{1}. Choose $r>0$ so that
$C_0:=\sup\{\rho(x):\|x\|\le r\}<\infty.$ Then
$$\textstyle
\E\sum\limits_{k=1}^n|\Delta^2\varphi(f_{k-1},df_k)|\le \E\rho(f_n)\le C_0
\qquad\text{whenever
$\|f_n\|_\infty\le r,$}
$$
where $\|g\|_\infty=\max_{\eta\in \Gamma^n}|g(\eta)|$ as usual.
Hence, for an arbitrary Walsh-Paley martingale $(f_0,\ldots,f_n)$, $p$-homogeneity implies that
\begin{equation}\label{ellinfty}
\textstyle
\E\sum\limits_{k=1}^n|\Delta^2\varphi(f_{k-1},df_k)| \le C_1\|f_n\|_\infty^p
\end{equation}
where $C_1=C_0/r^p$.

Now, fix any $X$-valued Walsh-Paley martingale $(f_0,\ldots,f_n)$ with
$\E \|f_n\|^p=1$. (By $p$-homogeneity, it suffices to prove \eqref{E:2} for such martingales.)
Let $\eta\in\Gamma^n$. We define
$m_0(\eta)=0$ and, for any integer $r\ge 1$,
$$
M_r(\eta)=\bigl\{0\le k<n: \max\{\|f_k(\eta)\pm df_{k+1}(\eta)\|\}>2^r\bigr\}
$$
and
\[
m_r(\eta)=\begin{cases}
\min M_r(\eta) &\text{if $M_r(\eta)\ne\emptyset$,}\\
n &\text{if $M_r(\eta)=\emptyset$.}
\end{cases}
\]
For each $m\in\{0,\ldots,n\}$, the set
$\{m_r=m\}$ belongs to $\Sigma_m$ by
Remark~\ref{O:mart}(f). (Thus the functions $m_r$ are so-called ``stopping times''.)
Hence it can be written in the form
$$
\{m_r=m\}=A_{r,m}\times\Gamma^{n-m},
\qquad\text{where $A_{r,m}\subset\Gamma^m$.}
$$
We have
\begin{align*}
\E \sum_{m_{r-1}<k\le m_r} & |\Delta^2\phi(f_{k-1},df_k)| \\
&=
\sum_{m=0}^{n-1}\, \int\limits_{\{m_{r-1}=m\}} \sum_{m<k\le m_r(\eta)}
|\Delta^2\phi(f_{k-1}(\eta),df_k(\eta))|\,d\p_n(\eta) \\
&= \sum_{m=0}^{n-1}\, \int\limits_{A_{r-1,m}}\left(
\,\,\int\limits_{\Gamma^{n-m}}
|\Delta^2\phi(f_{k-1}(\omega,\xi),df_k(\omega,\xi))|\,d\p_{n-m}(\xi)
\right)\,d\p_m(\omega)\,.
\end{align*}
The expression in parentheses can be seen as
\begin{equation}\label{exp}
\textstyle
\E \sum\limits_{m<k\le n} |\Delta^2\phi(g_{k-1}(\omega,\cdot),dg_k(\omega,\cdot))|\,,
\end{equation}
where
$$
g_k(\omega,\xi)=\begin{cases}
f_k(\omega,\xi) &\text{if $m<k\le m_r(\omega,\xi)$,} \\
f_{m_r(\omega,\xi)}(\omega,\xi) &\text{if $m_r(\omega,\xi)<k\le n$.}
\end{cases}
$$
Since $(g_k(\omega,\cdot))_{k=m}^n$ is a Walsh-Paley martingale by
Remark~\ref{O:mart}(d), and the definition of $m_r$ implies
$\|f_{m_r(\eta)}(\eta)\|=\|f_{m_r(\eta)-1}(\eta)+df_{m_r(\eta)}(\eta)\|\le 2^r$,
we can majorize the expression \eqref{exp} (using \eqref{ellinfty}) by
\[
C_1\|g_n(\omega,\cdot)\|_\infty^p=
C_1\|f_{m_r(\omega,\cdot)}(\omega,\cdot)\|_\infty^p\le C_1 2^{rp}.
\]
Thus
\begin{align*}\textstyle
\E\sum\limits_{k=1}^n|\Delta^2\varphi(f_{k-1},df_k)|
&\textstyle =
\sum\limits_{r=1}^\infty \,\E \sum\limits_{m_{r-1}<k\le m_r}  |\Delta^2\phi(f_{k-1},df_k)|
\\
&\textstyle
\le C_1\sum\limits_{r=1}^\infty 2^{rp} \sum\limits_{m=0}^{n-1} \p_m(A_{r-1,m}) =
C_1 \sum\limits_{r=1}^\infty 2^{rp} \,\p(m_{r-1}<n).
\end{align*}
Now, for $r>1$, Remark~\ref{O:mart}(f) implies that
$$
\p(m_{r-1}<n) \le
\p\bigl(\,\max\limits_{1\le k<n}\max\{\|f_k\pm df_{k+1}\|\}>2^{r-1}\bigr)
\le 2 \p\bigl(\,\max\limits_{0\le k\le n}\|f_k\|>2^{r-1}\bigr).
$$
This gives (via Observation~\ref{obs} and Lemma~\ref{L:mart}(d))
\begin{align*}
\textstyle
\E\sum\limits_{k=1}^n|\Delta^2\varphi(f_{k-1},df_k)|
&\textstyle
\le C_1 2^p + C_1 2\sum\limits_{r=2}^\infty 2^{rp}\,
\p(\max\limits_{0\le k\le n}\|f_k\|>2^{r-1}) \\
&\textstyle
= C_1 2^p + C_1 2^{p+1}\sum\limits_{j=1}^\infty 2^{jp}\,
\p(\max\limits_{0\le k\le n}\|f_k\|>2^j) \\
&\textstyle\le C_1 2^p + C_1\frac{2^{2p+1}}{2^p-1}\,
\E \max\limits_{0\le k\le n}\|f_k\|^p \\
&\le C_2 \bigl(1+\E\|f_n\|^p\bigr)=2C_2\,,
\end{align*}
where $C_2$ is a suitable constant depending only on $p$. Thus
\eqref{E:2} holds.

The converse follows trivially from Lemma~\ref{1} by putting
$\rho(x)=C\|x\|^p$.
\end{proof}

\begin{corollary}\label{p-homog}
Suppose $p\ge1$. Then every positively $p$-homogeneous delta-convex
function $\phi\colon X\to \R$ has a control function which is positively
$p$-homogeneous.
\end{corollary}

\begin{proof}
The case $p=1$ was proved in \cite[Lemma~1.21]{Ve-Za}. Assume $p>1$. By Lemma~\ref{2},
\eqref{E:1} holds with $\rho(x)=C\|x\|^p$. By the proof of
Lemma~\ref{1}, the formula \eqref{control} defines a positively $p$-homogeneous
control function for $\phi$.
\end{proof}

\begin{remark}
Let us remark that natural analogues of Lemma~\ref{1}, Lemma~\ref{2} and
Corollary~\ref{p-homog} hold also for mappings $\Phi\colon X\to Y$
(instead of functions $\phi\colon X\to\R$), where
``delta-convex function'' is replaced by
``delta-convex
mapping'' (as defined in \cite{Ve-Za}) and, in the terms involving $\Phi$,
the absolute value is replaced by the
norm of $Y$. This follows from
\cite[Proposition~1.13]{Ve-Za}.
\end{remark}

\begin{definition}
Let $X$ and $Y$ be Banach spaces. We say that a linear operator
$T\colon X\to Y$ is a
{\em UMD-operator} if there exists a constant $C>0$ such that
\begin{equation}\label{UMD}
\textstyle
\E \|\sum_{k=1}^n \epsilon_k\, Tdf_k\|^2 \le C\, \E \|f_n\|^2
\end{equation}
whenever
$(f_0,\ldots,f_n)$ is an $X$-valued Walsh-Paley martingale and
$\epsilon_1,\ldots,\epsilon_n$ are numbers in $\{-1,1\}$.
We say that $X$ is a {\em UMD-space}
if the identity $I\colon X\to X$
is a UMD-operator.
\end{definition}

\begin{remark}\label{R:UMD}
\
\begin{enumerate}
\item[(a)] It is easy to see that a composition of two bounded linear
operators is a UMD-operator whenever at least one of them is. In particular, 
if at least one of $X,Y$ is a UMD-space, 
then each bounded linear
operator $T\colon X\to Y$ is a UMD-operator.
\item[(b)] Suppose $p>1$ is a real number. Then $X$ is a UMD-space 
if and only if
is a constant $c_p>0$ such that
$$
\textstyle
\E \|\sum_{k=1}^n \epsilon_k \,df_k\|^p \le c_p\, \E \|f_n\|^p
$$
whenever $\epsilon_k=\pm1$ ($1\le k\le n$) and $(f_0,\ldots,f_n)$ is a
 Walsh-Paley martingale.  Moreover, in this case the above inequality holds also for general (i.e. not necessarily Walsh-Paley) martingales.
(See p.67 and Lemma~7.1 in \cite{B-LNM}.)
\item[(c)] Every UMD-space is superreflexive.
(See e.g.\ \cite[p.222]{Pi-LNM} or \cite[Proposition~2]{aldous}.)
\end{enumerate}
\end{remark}

For us the following result, which was proved in
\cite{wenzel}, will be important. Let us remark that a similar result for general martingales in UMD-spaces was proved
by Burkholder~\cite{B} and, as remarked in
\cite[p.502]{BD}, his proof can be easily modified to prove
the same for UMD-operators defined using general martingales.

\begin{fact}\label{UMD-MT}
Let $T\colon X\to Y$ be a UMD-operator between Banach spaces $X,Y$. Then
there exists a constant $C>0$ such that \eqref{UMD} holds whenever
$(f_0,\ldots,f_n)$ is an $X$-valued Walsh-Paley martingale and
$(\epsilon_1,\ldots,\epsilon_n)$ is a predictable sequence of
$\{\pm1\}$-valued functions
(i.e., each $\epsilon_k$ is $\Sigma_{k-1}$-measurable).
\rm (See \cite{wenzel}.)
\end{fact}

\begin{theorem}\label{3}
Let $q$ be a continuous quadratic form on $X$ and $T\colon X\to X^*$ the
symmetric operator that generates $q$. Then $q$ is delta-convex if and
only if $T$ is a UMD-operator.
\end{theorem}

\begin{proof}
Let $T$ be a UMD-operator, let $(f_0,\ldots,f_n)$ be an $X$-valued Walsh-Paley martingale.
Then, for $\eta\in\Gamma^n$ and $1\le k\le n$, we have
$$
|\Delta^2 q(f_{k-1}(\eta),df_k(\eta))|
= 2|q(df_k(\eta))|
= 2\epsilon_k(\eta) q(df_k(\eta))
=2\epsilon_k(\eta) \langle df_k(\eta),Tdf_k(\eta)\rangle,
$$
where $\epsilon_k(\eta)=1$ if $q(df_k(\eta))\ge0$,
$\epsilon_k(\eta)=-1$ if $q(df_k(\eta))<0$.
Observe that $\epsilon_k$ is $\Sigma_{k-1}$-measurable by
Remark~\ref{O:mart}(e) since $q$ is an even function; in other words,
the sequence $(\epsilon_1,\ldots,\epsilon_n)$ is predictable. Using Lemma~\ref{L:mart}, we can
write
\begin{align*}
\textstyle
\E \sum\limits_{k=1}^n |\Delta^2 q(f_{k-1},df_k)|
&\textstyle=
2 \E \sum\limits_{k=1}^n  \langle df_k,\epsilon_k Tdf_k\rangle
=
2\E \langle\sum\limits_{j=1}^n df_j\,,\sum\limits_{k=1}^n \epsilon_k Tdf_k\rangle\\
&\textstyle=
2\E \langle f_n-f_0,\sum\limits_{k=1}^n \epsilon_k Tdf_k\rangle
=
2\E \langle f_n,\sum\limits_{k=1}^n \epsilon_k Tdf_k\rangle\\
&\textstyle\le
2 (\E\|f_n\|^2)^{1/2} \bigl(\E\|\sum\limits_{k=1}^n \epsilon_k Tdf_k\|^2\bigr)^{1/2}
\le
2\sqrt{C}\,\, \E\|f_n\|^2\,,
\end{align*}
where $C$ is the constant from Fact~\ref{UMD-MT}.

For the converse, suppose that $q$ is delta-convex. Consider any
$X$-valued Walsh-Paley martingale $(f_0,\ldots,f_n)$ with
$\E\|f_n\|^2=1$, and numbers $\epsilon_k\in\{-1,1\}$ ($1\le k\le n$). 
It is easy to see that
there exists $h_n\colon\Gamma^n\to X$ such that
$\E\|h_n\|^2=1$ and
$$\textstyle
\bigl(\|\sum\limits_{k=1}^n \epsilon_k Tdf_k\|^2\bigr)^{1/2}\le
2 \E\langle h_n, \sum\limits_{k=1}^n \epsilon_k Tdf_k\rangle\,;
$$
indeed, 
denoting $g=\sum_{k=1}^n \epsilon_k Tdf_k$,
one can put
$h(\eta):=t(\eta) v(\eta)$,
where 
$t\colon\Gamma^n\to[0,\infty)$ satisfies
$\E t^2=1$ and $(\E\|g\|^2)^{1/2}=\E(t\cdot\|g\|)$, and
$v(\eta)\in S_X$ is such that $\|g(\eta)\|\le 2\langle
v(\eta),g(\eta)\rangle$.
Let $(h_0,\ldots,h_n)$ be the Walsh-Paley martingale given by $h_n$
(i.e., $h_k=\E(h_n|\Sigma_k)$, $1\le k\le n$). Then
Lemma~\ref{2} and
the identities 
$$\langle x,Ty\rangle = (1/4)\bigl(q(x+y)-q(x-y)\bigr)\ \ \ 
\text{and}\ \ \ q(y)=(1/2)\Delta^2q(x,y)$$ 
imply
\begin{align*}
\textstyle
\bigl(\|\sum\limits_{k=1}^n \epsilon_k Tdf_k\|^2\bigr)^{1/2}&
\textstyle\le
2 \E\langle h_n, \sum\limits_{k=1}^n \epsilon_k Tdf_k\rangle\,=
\qquad\text{(as above)}\\
&\textstyle=
2\E \sum\limits_{k=1}^n \langle dh_k, \epsilon_k Tdf_k\rangle\le
2\E \sum\limits_{k=1}^n \bigl| \langle dh_k, Tdf_k\rangle
\bigr|\\
&\textstyle\le
\frac12 \E \sum\limits_{k=1}^n \bigl|q(d(f_k+h_k))\bigr| +
\frac12 \E \sum\limits_{k=1}^n \bigl|q(d(f_k-h_k))\bigr| \\
&\textstyle\le
\frac14 C\, \E\|f_n+h_n\|^2 + \frac14 C\, \E \|f_n-h_n\|^2 \\
&\textstyle\le
\frac12 C \E\bigl(\|f_n\|+\|h_n\|\bigr)^2 
\le C\bigl(\E\|f_n\|^2 +\E\|h_n\|^2\bigr)\le 2C.
\end{align*}
Thus $T$ is a UMD-operator.
\end{proof}

The following theorem is an immediate consequence of Theorem~\ref{3} and
Remark~\ref{R:UMD}(a).

\begin{theorem}\label{4}  Let $X$ be a Banach space.  Then every continuous quadratic form on $X$ is
delta-convex if and only if every symmetric operator $T:X\to X^*$ is a UMD-operator.
In particular, if $X$ is a UMD-space, then every continuous quadratic form on $X$ is
delta-convex.
\end{theorem}

\begin{theorem}\label{T:counterexample}
There exists a continuous quadratic form on $\ell_1$ which is not
delta-convex.
\end{theorem}

\begin{proof}
Let $J\colon \ell_1\to\ell_\infty$ be an isometric embedding (recall that every separable Banach space isometrically embeds into $\ell_{\infty}$). Consider
the continuous quadratic form on $\ell_1=\ell_1\oplus_1\ell_1$, given by
$$q(x,y)=\langle y, Jx \rangle +  \langle x, Jy \rangle\,,
\qquad x,y\in\ell_1\,,$$
which is generated by the symmetric operator $T(x,y)=(Jy,Jx)$. If $q$ were
delta-convex, $T$ would be a UMD-operator. But then $J$
would be a UMD-operator; consequently, $\ell_1$ would be a UMD-space.
But this is false by Remark~\ref{R:UMD}(c).
\end{proof}

Let us conclude with a simple but useful proposition.

\begin{proposition}\label{P:d.c.} Let $p>0$.
Let $\phi\colon X\to\R$ be a $p$-homogeneous function on a Banach space
$X$. Then $\phi$ is delta-convex if and only if $\phi$ is delta-convex
on a convex neighborhood of the origin.
\end{proposition}

\begin{proof}
Let $\phi$ be delta-convex on a convex neighborhood $U$ of the origin, and let 
$\psi\colon U\to\R$ be a corresponding control function.
There exists $\delta>0$ such that
$\psi$ is bounded on $\delta B_X$.
A simple homogeneity argument shows that $\phi$ is delta-convex on each
$r B_X$ ($r>0$)
with a bounded control function of the form
$\rho(x)=c_1\psi(c_2 x)$.
Then $\phi$ is delta-convex on $X$
by \cite[Theorem~16]{Kop-Maly}.
\end{proof}

%%%%%%%%%%%%%%%%%%%%%%%%%%%%%%%%%%%%%%%%%%%%%%%%%%%%%%%%%%%%%%%

\section{Further results and open problems}\label{S:three}

We shall consider the following three properties of a Banach space $X$, defined already
in Introduction.

\begin{enumerate}
\item[(D)] {\em Each continuous quadratic form on $X$ is
          delta-semidefinite.}

\item[(dc)] {\em Each continuous quadratic form on $X$ is
          delta-convex.}

\item[(Cdc)] {\em Each $C^{1,1}$ function $f\colon X\to\R$ is
          delta-convex.}
\end{enumerate}

\medskip\noindent
Recall that a function (or mapping) $f$ is
$C^{1,1}$ if the Fr\'echet
derivative $f'(x)$ exits for each $x$ and the mapping
$f'$ is Lipschitz.

We have seen that (D) passes to quotients (Corollary~\ref{C:quotient}).  Let us observe the
same result for properties (dc) and (Cdc).

\begin{lemma}\label{dcquotient} 
If $X$ is a Banach space with property (dc) 
(respectively, (Cdc)), then for any closed subspace $E$ of $X$, the quotient 
$X/E$ has property (dc) (respectively, (Cdc)).
\end{lemma}

\begin{proof} Let $Q\colon X\to X/E$ be the quotient map and let $f\colon X/E\to\mathbb R$ 
be a continuous function such that $f\circ Q$ is
delta-convex.  We show that $f$ is delta-convex (which proves both assertions).  
Let $\psi\colon X\to\mathbb R$ be a continuous function such that $\psi\pm f$ is convex; 
we can assume $\psi\ge 0$.  Define $\hat\psi\colon X/E\to\mathbb R$ by
$\hat\psi(y)=\inf\{\psi(x): Qx=y\}.$  Then it is easy to prove that
$\hat\psi\pm f$ is convex. Moreover, the (convex)
function $\hat\psi$ is continuous since it is easily seen to be
bounded on a neighborhood of the origin.
\end{proof}

Since each continuous quadratic form is $C^{1,1}$
(by~Fact~\ref{F:p'}),  we always have the
implications
\[
({\rm D})\ \Longrightarrow\ ({\rm dc})\
\Longleftarrow\ ({\rm Cdc})\,.
\]
As we shall see in the next theorem, no two of the above three
properties are equivalent.

Let us start with the following corollary of 
\cite[Theorem~11]{Du-Ve-Za}. A norm on $X$ is said to have modulus of
convexity of power type 2 if,
for some $c>0$, $\delta_{X}(\epsilon)\ge c\cdot\epsilon^2$
whenever $\epsilon\in(0,2]$ (where $\delta_X$ is the usual modulus of convexity;
see e.g.\ \cite{Lin-Tz}). 

\begin{fact}\label{F:powertype}
Let $X$ be a Banach space that admits a uniformly convex renorming 
with modulus of convexity of power type 2. Then $X$ satisfies (Cdc) and
hence also (dc).
\end{fact}

By an $L_p(\mu)$ space we mean an infinite-dimensional space $L_p(\Omega,\Sigma,\mu)$ where
$(\Omega,\Sigma,\mu)$ is a positive measure space. This class includes the spaces $L_p(0,1)$ and $\ell_p$. For such spaces,
we have the following theorem which summarizes
results of \cite{Du-Ve-Za}, \cite{Zel} and of the present paper.

\begin{theorem}\label{T:L_p}
Let $X$ be an infinite-dimensional $L_p(\mu)$ space with
$1\le p \le \infty$.
\begin{enumerate}
\item[(a)]
$X$ satisfies {\rm(D)} if and only if $p\ge2$.
\item[(b)]
$X$ satisfies {\rm(Cdc)} if and only if $1<p\le2$.
\item[(c)]
$X$ satisfies {\rm(dc)} if and only if $p>1$.
\end{enumerate}
\end{theorem}

\begin{proof}
{\em(a)} If $p\ge2$ then $L_p(\mu)$ satisfies (D) by
Theorem~\ref{T:D}(f). If $p<2$ then $L_p(\mu)$ fails (D) since it
contains a complemented copy of $\ell_p$ which fails (D) by
Corollary~\ref{C:badp}.

{\em(b)} If $1<p\le2$, then the standard norm on $X=L_p(\mu)$ has modulus of
convexity of power type 2 (see\ \cite[Corollary~V.1.2]{DGZbook}).
By Fact~\ref{F:powertype}, each such space satisfies (Cdc).
$L_1(\mu)$ fails (Cdc) since it fails (dc) (see (c) below). Now, let $2<p<\infty$.
For such $p$, M.~Zelen\'y \cite{Zel} proved that $\ell_p$ fails (Cdc);
thus $L_p(\mu)$ fails (Cdc), too. Finally, to see that also $L_\infty(\mu)$ fails
(Cdc), it suffices to show that $L_\infty(0,1)$ fails (Cdc); indeed, the
spaces $L_\infty(0,1)$ and $\ell_\infty$ are isomorphic by \cite{pelcz},
and $L_\infty(\mu)$ contains a complemeted copy of $\ell_\infty$.
By \cite[Corollary~2.f.5]{Lin-Tz}, $(\ell_4)^*=\ell_{4/3}$ isometrically
embeds in $L_1(0,1)$; consequently, $\ell_4$ is isomorphic to a quotient
of $L_\infty(0,1)$. Hence $L_\infty(0,1)$ fails (Cdc) by
Lemma~\ref{dcquotient}, since we already know that $\ell_4$ fails (Cdc).

{\em(c)} $L_1(\mu)$ fails (dc) since it contains a complemented copy of $\ell_1$
which fails (dc) by Theorem~\ref{T:counterexample}. For $p>1$,
the space $L_p(\mu)$ satisfies (dc) since, by (a) and (b), it satisfies
(D) (if $p\ge2$) or (Cdc) (if $p\le2$).  Alternately one may observe that 
$L_p(\mu)$ is a UMD-space if $1<p<\infty$ using Remark \ref{R:UMD}(b) and 
then apply Theorem~\ref{4}.
\end{proof}

\begin{remark}
By Theorem~\ref{T:L_p},
${\rm(dc)}\not\Rightarrow{\rm(D)}$ and also
${\rm(dc)}\not\Rightarrow{\rm(Cdc)}$.
It is interesting to compare the second non-implication with the following
result from \cite{Du-Ve-Za} about vector-valued mappings:
{\em
every Banach space-valued continuous quadratic mapping on $X$ is delta-convex
if and only if every Banach space-valued $C^{1,1}$ mapping on $X$ is
delta-convex.}
\end{remark}

\subsection*{Delta-convex functions via delta-convex curves}
There is another corollary to the above results. It is connected with
Problem~6 in \cite{Ve-Za}.
In that paper,
delta-convex mappings between Banach spaces
(a generalization of delta-convex functions)
were defined and widely studied. We do not state
the definition here; it can be found also in \cite{Du-Ve-Za} together with a survey of
principal results. We confine ourselves to
stating an equivalent definition (see~\cite[Theorem~2.3]{Ve-Za})
of a delta-convex mapping of one real
variable.

\begin{definition}
Let $I\subset\R$ be an open interval, $X$ be a Banach space. A mapping
$\phi\colon I\to X$ is {\em delta-convex} on $I$ if the right derivative $\phi'_+(t)$
exists at each $t\in I$ and the mapping $\phi'_+$ has bounded variation
on each compact subinterval of $I$.
\end{definition}

For real-valued functions, Problem~6 in \cite{Ve-Za} asks: {\em suppose that
$X$ is a Banach space and $f\colon X\to\R$ is a function such that
$f\circ\phi$ is a delta-convex function on $(0,1)$
whenever $\phi\colon(0,1)\to X$ is a delta-convex mapping;
is then $f$ locally delta-convex?} The following example answers in
negative this problem. (Let us remark that the
vector-valued case was solved in negative already in
\cite{Du-Ve-Za}.)

\begin{example}\label{E:composition}
Let $X$ be an infinite-dimensional $L_p(\mu)$ space where either $p=1$
or $2<p<\infty$. Then there exists a continuous function $f\colon X\to\R$ such that
$f$ is delta-convex on no neighborhood of $0$, and
$f\circ\phi$ is  delta-convex on $(0,1)$ for each delta-convex mapping
$\phi\colon(0,1)\to X$.
\end{example}

\begin{proof}
The case $p=1$. By Theorem~\ref{T:L_p}(b), there exists a continuous
quadratic form $q$ on $X$ such that $q$ is not delta-convex. By
Proposition~\ref{P:d.c.}, $q$ is delta-convex on no neighborhood of $0$.
By Proposition~14 in \cite{Du-Ve-Za}, $q\circ\phi$ is delta-convex for
each delta-convex ``curve'' $\phi$. Thus we can put $f=q$.

The case $2<p<\infty$ follows in a similar way using \cite{Zel} instead of
Theorem~\ref{T:L_p}. Indeed, by
\cite{Zel}, there exists a $C^{1,1}$ function
$g\colon\ell_p\to\R$ that is not delta-convex. A careful look at the
proof in \cite{Zel} shows that the function constructed therein
is d.c.\ on no neighborhood of $0$.
Consider $\ell_p$ as a complemented subspace of $X=L_p(\mu)$
and extend $g$ to a $C^{1,1}$ function on the whole $X$ by $f=g\circ P$
where $P$ is a bounded linear projection of $X$ onto $\ell_p$. Then $f$
has the desired property by \cite[Proposition~14]{Du-Ve-Za} again.
\end{proof}

It is natural to ask the following

\begin{problem}
Does there exist a function $f$ as in Example~\ref{E:composition} for
each at least two-dimensional Banach space $X$? Or, at least, for each
infinite-dimensional Banach space $X$?
\end{problem}

\subsection*{Stability with respect to direct sums}
Consider two Banach spaces $X_1$ and $X_2$. Since
$(X_1\oplus X_2)^*=X_1^*\oplus X_2^*$, each bounded linear operator
operator $T\colon X_1\oplus X_2 \to (X_1\oplus X_2)^*$ can be
represented as an operator-valued matrix
$T=\textstyle\left( \begin{array}{cc}
 T_{11} &  T_{12} \\
 T_{21}  &  T_{22}
\end{array}\right)$ where
$T_{ij}\colon X_j\to X_i^*$.
It is an easy exercise to verify that:
\begin{enumerate}
\item[(I)] $T$ is factorizable through a Hilbert space if and only if each
$T_{ij}$ is;
\item[(II)] $T$ is a UMD-operator if and only if each $T_{ij}$ is;
\item[(III)] $T$ is symmetric if and only if
$T_{ij}^*=T_{ji}$ on $X_i$ whenever $i,j\in\{1,2\}$ (equivalently, $T_{11}$ and $T_{22}$ 
are symmetric and $T^*_{12}=T_{21}$ on $X_1$).
\end{enumerate}
Hence we have the following consequence of Corollary~\ref{C:quotient} and
Theorem~\ref {3}.
(Given Banach spaces $X$ and $Y$, we denote by $\mathcal{L}(X,Y)$ the set of
all bounded linear operators from $X$ into $Y$.)

\begin{corollary}\label{C:sum} Let $X_1$ and $X_2$ be Banach spaces.
Then:
\begin{enumerate}
\item[(a)] $X_1\oplus X_2$ has the property (D) if and only if $X_1$
and $X_2$ have (D) and every element of
$\mathcal{L}(X_1,X_2^*)$ is factorizable through
a Hilbert space; 
\item[(b)] $X_1\oplus X_2$ has the property (dc) if and only if $X_1$
and $X_2$ have (dc) and every element of
$\mathcal{L}(X_1,X_2^*)$ is a UMD-operator.
\end{enumerate}
In particular, if $X$ is isomorphic to $X^2$ then
\begin{enumerate}
\item[(a')] $X$ has the property (D) if and only if every element of
$\mathcal{L}(X,X^*)$ is factorizable through
a Hilbert space; 
\item[(b')] $X$ has the property (dc) if and only if every element of
$\mathcal{L}(X,X^*)$ is a UMD-operator.
\end{enumerate}

\end{corollary}

The following Corollary is immediate using Remark \ref{R:UMD}(a) for part (b).

\begin{corollary}
Let $X$ be a Banach space.
\begin{enumerate}
\item[(a)] If $X$ satisfies (D) and $H$ is a Hilbert space, 
then $X\oplus H$ satisfies (D).
\item[(b)] If $X$ satisfies (dc) and $U$ is a UMD-space, 
then $X\oplus U$ satisfies (dc).
\end{enumerate}
\end{corollary}

\begin{observation}\label{adjoint}
The adjoint $T^*$ is a UMD-operator if and only if $T$ is a UMD-operator. In particular, $X$
is a UMD-space if and only if $X^*$ is. \\
\rm
To see this, note that $T\colon X\to Y$ is a UMD-operator if and only if
the operators
$$\textstyle
T_{n,\epsilon}:=\sum_{k=1}^n \epsilon_k T(E_{k,X}-E_{k-1,X})
\colon L_2(\Gamma^n,X)\to L_2(\Gamma^n,Y),\quad\;
n\in\N,\;\epsilon\in\{-1,1\}^n,
$$
are equi-bounded, where $E_{k,X}:=\E(\cdot|\Sigma_k)\colon
L_2(\Gamma^n, X)\to L_2(\Gamma^n, X)$. But then also the corresponding adjoints are
equi-bounded. Moreover, it is easy to see that
$$\textstyle
T_{n,\epsilon}^*=
\sum_{k=1}^n \epsilon_k (E_{k,X^*}-E_{k-1,X^*}) T^*=
\sum_{k=1}^n \epsilon_k T^*(E_{k,Y^*}-E_{k-1,Y^*})$$ which means that
$T^*$ is a UMD-operator. The reverse implication follows easily from this one.
\end{observation}

\begin{proposition}\label{XxX*}
Let $X$ be a Banach space. Let $Q$ be the continuous quadratic
form on $X\oplus X^*$, given by $Q(x,x^*)=x^*(x)$. Then the following three assertions are
equivalent:
\begin{enumerate}
\item[(i)] $X\oplus X^*$ satisfies (D);
\item[(ii)] $Q$ is delta-semidefinite;
\item[(iii)] $X$ is isomorphic to a Hilbert space.
\end{enumerate}
And also the following three assertions are equivalent:
\begin{enumerate}
\item[(i')] $X\oplus X^*$ satisfies (dc);
\item[(ii')] $Q$ is delta-convex;
\item[(iii')] $X$ is a UMD-space.
\end{enumerate}
\end{proposition}

\begin{proof}
The implications (iii)$\Rightarrow$(i)$\Rightarrow$(ii) are obvious.
Assume (ii). Since $Q$ is generated by the symmetric operator
$T\colon X\oplus X^*\to X^*\oplus X^{**}$, $T(x,x^*)=\frac12(x^*,x)$, it
follows from Theorem~\ref{T:kalton} that the identity $I\colon X\to X$ factors
through a Hilbert space. It is easy to see that this implies (iii). (Indeed,
if $I=BA$ is a factorization through a Hilbert space $H$, then $AB$ is a
bounded linear projection onto a closed subspace $H_0=A(X)$ of $H$. Then $A$ is
a linear isomorphism between $X$ and $H_0$.)

The implication (i')$\Rightarrow$(ii') is obvious. If (ii') holds, then (as
above, via Theorem~\ref{3}) the identity $I\colon X\to X$ is a UMD-operator, 
which gives
(iii'). Finally, if (iii') holds, then $X\oplus X^*$ is a UMD-space
(indeed, it suffices to apply (II) before
Corollary~\ref{C:sum} to the identity operator of $X\oplus X^*$, taking
into account Observation~\ref{adjoint}). Hence (i') holds by Theorem~\ref{4}.
\end{proof}

As far as we know, the following question is open.

\begin{problem}\label{DDD}
Is the property (D) stable with respect to making direct sums of two
spaces? Equivalently, if Banach spaces $X_1$ and $X_2$ have property
(D), does it imply that each $S\in\mathcal{L}(X_1,X_2^*)$ is factorizable
through a Hilbert space?
\end{problem}

We conjecture that the answer is negative, but we do not know any
counterexample. However, the following observation shows that a possible
counterexample cannot be found by using only spaces provided by
Theorem~\ref{T:D}.

\begin{observation}
Let each of given two Banach spaces $X_1$ and $X_2$ satisfy at least one
of the conditions (a)--(i) in Theorem~\ref{T:D}. Then $X_1\oplus X_2$ has
(D).
\end{observation}

\begin{proof}
By Corollary~\ref{C:sum}, it suffices to show that every operator
$S\in\mathcal{L}(X_1,X_2^*)$ factors through a Hilbert space. By the
proof of Theorem~\ref{T:D}, each of the spaces $X_i$ (i=1,2) has at
least one of the following three properties:
\begin{enumerate}
\item[($\alpha$)] $X_i$ has type 2;
\item[($\beta$)] $X_i^*$ has cotype 2, and $X_i$ has the approximation
property;
\item[($\gamma$)] $X_i^*$ has cotype 2, and $X_i$ is a Banach lattice.
\end{enumerate}
First observe that this implies that $X_2^*$ has cotype 2 (see e.g.\
\cite[Proposition~3.2]{Pisier}). Now, if $X_1$ satisfies ($\alpha$),
apply \cite[Corollary~3.6]{Pisier}; if $X_1$ satisfies ($\beta$), apply
\cite[Theorem~4.1]{Pisier}; if $X_1$ satisfies ($\gamma$), apply
\cite[Theorems~8.17 and 8.11]{Pisier}.
\end{proof}

Note that Proposition~\ref{XxX*} implies that Problem~\ref{DDD} will have
a negative answer once the following problem is solved in negative.

\begin{problem}
Let $X$ and $X^*$ satisfy (D). Does it imply that $X$ is isomorphic
to a Hilbert space?
\end{problem}

For the property (dc) we have the following theorem.

\begin{theorem}
The property (dc) is not stable under making direct sums.
\end{theorem}

\begin{proof} 
By \cite{bourgain}, there exists a Banach lattice $X$ 
such that $X$ is not a UMD-space and $X$ satisfies an upper-3 estimate
and a lower-4 estimate (see \cite{Lin-Tz} for definitions). By 
\cite[Theorem~1.f.7]{Lin-Tz}, $X$ is 2-convex and 5-concave. Then $X$
admits a uniformly smooth renorming with modulus of smoothness of 
power type 2 by \cite[Theorem~1.f.1]{Lin-Tz}, which implies that
$X$ has (D) (see Theorem~\ref{T:D}(h)), and hence (cd). By duality (see
\cite[p.63]{Lin-Tz}), $X^*$ admits a uniformly convex renorming with
modulus of convexity of power type 2; consequently, $X^*$ has (dc) by
 Fact~\ref{F:powertype}. By Proposition~\ref{XxX*}, $X\oplus X^*$ fails (dc) since $X$ is not 
a UMD-space.
\end{proof}

%%%%%%%%%%%%%%%%%%%%%%%%%%%%%%%%%%%%%%%%%%%%%%%%%%%%%%%%%%%%%%%

%\subsection*{Acknowledgment}

%%%%%%%%%%%%%%%%%%%%%%%%%%%%%%%%%%%%%%%%%%%%%%%%%%%%%%%%%%%%%%%

\end{document}